\documentclass[preprint,12pt]{elsarticle}

\usepackage{amssymb}
\usepackage{amsthm}

\usepackage{amsmath}

\newtheorem{theorem}{Theorem}[section]
\newtheorem{lemma}[theorem]{Lemma}
\newtheorem {proposition} [theorem] {Proposition}

\theoremstyle{definition}

\newtheorem{example}[theorem]{Example}

\theoremstyle{remark}
\newtheorem{remark}[theorem]{Remark}

\numberwithin{equation}{section}

\journal{XXXXX}

\begin{document}

\begin{frontmatter}



\title{Conditions to the existence of center in planar systems and center for Abel equations}


\author[a]{Anderson L. A. de Araujo}
\ead{anderson.araujo@ufv.br}
\author[a]{Ab\'ilio Lemos}
\ead{abiliolemos@ufv.br}
\author[a]{Alexandre M. Alves}
\ead{amalves@ufv.br}

\address[a]{Departamento de Matem\'atica, Universidade Federal de Vi\c cosa, 36570-900, Vi\c cosa (MG), Brazil}

\begin{abstract}
Abel equations of the form $x'(t)=f(t)x^3(t)+g(t)x^2(t)$, $t \in [-a,a]$, where $a>0$ is a constant, $f$ and $g$ are continuous functions, are of interest because of their close relation to planar vector fields. If $f$ and $g$ are odd functions, we prove, in this paper, that the Abel equation has a center at the origin. We also consider a class of polynomial differential equations $\dot{x} =  -y+P_n(x,y)$ and  $\dot{y} =  x+Q_n(x,y)$, where $P_n$ and $Q_n$ are homogeneous polynomials of degree $n$. Using the results obtained for Abel's equation, we obtain a new subclass of systems having a center at the origin.
\end{abstract}


\end{frontmatter}


\section{Introduction and main results}

\subsection{Historical Aspect}

Let the planar system

\begin{equation}\label{S1}
\displaystyle{\begin{array}{ccc}
\dot{x} & = &  -y+P(x,y)\\
 \dot{y} & = &  x+Q(x,y),\\
\end{array}}
 \end{equation}
where $P(x,y)$ and $Q(x,y)$ are polynomials, without constant term, of maximum degree n. The singular point $(0,0)$ is a center, if surrounded by closed trajectories; or a focus, if surrounded by spirals. The classical center-focus problem consists in distinguishing when a  singular point is either a center or a focus. The problem started with Poincar\'e \cite{P} and Dulac \cite{Du}, and, in the present days, many questions remain open. The basic results were obtained by A. M.
Lyapunov \cite{Li}. He proved that if $P(x,y)$ and $Q(x,y)$ satisfy an infinite sequence of recursive conditions, then \eqref{S1} has a center to the origin. He also presented conditions for the origin of the system \eqref{S1} to be a focus.

If we write $P(x,y)=\sum_{i=1}^{l}P_{m_i}(x,y)$ and $Q(x,y)=\sum_{i=1}^{l}Q_{m_i}(x,y)$, where $P_{m_i}(x,y)$ and $Q_{m_i}(x,y)$ are homogeneous polynomials of degree $m_i\geq1$, then, from Hilbert's theorem on the finiteness of basis of polynomial ideals (\cite{Ka}, Theorem 87, p. 58), it follows that in the mentioned infinite sequence of recursive conditions only a finite number of conditions for center are essential. The others result from them.  

In this paper, we study a particular case of \eqref{S1}. Namely

\begin{equation}\label{eq.1}
\displaystyle{\begin{array}{ccc}
\dot{x} & = &  -y+P_n(x,y)\\
 \dot{y} & = &  x+Q_n(x,y),\\
\end{array}}
 \end{equation}
where $P_n(x,y)$ and $Q_n(x,y)$ are homogeneous polynomials of degree n.  

When $n = 2$, systems \eqref{eq.1} are quadratic polynomial differential systems (or simply
quadratic systems in what follows). Quadratic systems have been intensively studied over
the last 30 years, and more than a thousand papers on this issue have been published (see,
for example, the bibliographical survey of Reyn \cite{Re}).

%

A method for investigating if \eqref{eq.1} has a center at the origin is to transform the planar system into an Abel equation. 
In polar coordinates $(r,\theta)$ defined by $x = r\cos\theta, y = r\sin\theta$, the system (\ref{eq.1}) becomes
\begin{equation}\label{eq.2}
\displaystyle{\begin{array}{ccl}
\dot{r} & = &  A(\theta)r^n\\
 \dot{\theta} & = &1+B(\theta)r^{n-1},\\
\end{array}}
 \end{equation}
where 
\begin{equation}\label{eq.3}
\displaystyle{\begin{array}{ccl}
A(\theta) & = &  \cos\theta P_n(\cos\theta,\sin\theta) + \sin\theta Q_n(\cos\theta,\sin\theta),\\
B(\theta) & = & \cos\theta Q_n(\cos\theta,\sin\theta) - \sin\theta P_n(\cos\theta,\sin\theta).\\
\end{array}}
 \end{equation}
We remark that $A$ and $B$ are homogeneous polynomials of degree $n + 1$ in the variables
$\cos\theta$ and $\sin\theta$. In the region
$$R=\{(r,\theta):1+B(\theta)r^{n-l}>0\},$$
the differential system (\ref{eq.2}) is equivalent to the differential equation

\begin{equation}\label{eq.4}
\displaystyle{\frac{dr}{d\theta} =\frac{A(\theta)r^n}{1+B(\theta)r^{n-1}}}.
\end{equation}
It is known that the periodic orbits surrounding the origin of the system (\ref{eq.2}) does not intersect the curve $\theta = 0$ (see the Appendix of \cite{CL}). Therefore, these periodic orbits are contained in the region $R$. Consequently, they are also periodic orbits of equation (\ref{eq.4}).

The transformation $(r, \theta) \rightarrow (\gamma, \theta)$ with

\begin{equation}\label{eq.5}
\displaystyle{\gamma =\frac{r^{n-1}}{1+B(\theta)r^{n-1}}}
\end{equation}
is a diffeomorphism from the region R into its image. As far as we know, Cherkas was the first to use this transformation (see \cite{C}). If we write equation (\ref{eq.4}) in the variable $\gamma$, we obtain
\begin{equation}\label{eq.6}
\displaystyle{\frac{d\gamma}{d\theta} =-(n-1)A(\theta)B(\theta)\gamma^3 + [(n-1)A(\theta)-B'(\theta)]\gamma^{2}},
\end{equation}
which is a particular case of an Abel differential equation. We notice that $f(\theta)=-(n-1)A(\theta)B(\theta)$ and $g(\theta)=(n-1)A(\theta)-B'(\theta)$ are homogeneous trigonometric polynomials of
degree $2(n + 1)$ and $n + 1$, respectively.

Now the Center-Focus problem of planar system (\ref{eq.1}) has a translation in equation (\ref{eq.6}), that is, given $\gamma_0$ small enough, we look for necessary and sufficient conditions on $f(\theta)$ and $g(\theta)$ in order to assure that the solution of equation (\ref{eq.6}), with the initial condition $\gamma(0) = \gamma_0$, has the property that $\gamma(0) = \gamma(2\pi)$. We observe that this condition implies the periodicity of this solution, in particular, one has $\gamma(-\pi)=\gamma(\pi)$.

The equation $dx/dt = a(t)x^3 + b(t)x^2$ was studied in \cite{AGG}, where necessary and sufficient conditions were obtained for this equation has a center at the origin, but with $a(t)$ and $b(t)$, particular continuous functions (see Example \ref{ex3}). More results that ensure the existence a center at the origin for some subclasses of planar systems and for Abel equations were obtained in \cite{LLZ,LLLZ}.  

In this paper, some new results are obtained for planar systems \eqref{eq.1} and Abel equations \eqref{eq.6}.

\subsection{Results on planar systems}

\indent Consider the planar system of differential equations $(\ref{eq.1})$. We obtain sufficient conditions on $P_n(x,y)$ and $Q_n(x,y)$ that allow the system $(\ref{eq.1})$ to have a center at the origin.

\begin{theorem}\label{P1}
Let the planar system 
\begin{equation}\label{eq.13}
\displaystyle{\begin{array}{ccc}
\dot{x} & = &  -y+P_n(x,y)\\
 \dot{y} & = &  x+Q_n(x,y),\\
\end{array}}
 \end{equation} 
where $P_n(\cos\theta,\sin\theta)$ is an odd function and $Q_n(\cos\theta,\sin\theta)$ is an even  function in $C([0,2\pi])$. Then, the system (\ref{eq.13}) has a center at the origin.
\end{theorem}

\begin{example}\label{ex4}
Let the planar system 
\begin{equation}\label{eq.15}
\displaystyle{\begin{array}{ccc}
\dot{x} & = &  -y+ax^{N_1}y^{M_1}\\
 \dot{y} & = &  x+bx^{N_2}y^{M_2},\\
\end{array}}
 \end{equation} 
with $N_1+M_1=N_2+M_2$, $N_1,N_2$ be any nonnegative integers numbers, $M_1$ is an odd natural number, $M_2$ is an even natural number and $a,b$ are real numbers. Then, the planar system \eqref{eq.15} has a center at the origin. Indeed, it is easy to see that $P_n(\cos\theta,\sin\theta)=a\cos^{N_1}\theta\sin^{M_1}\theta$ is an odd function and $Q_n(\cos\theta,\sin\theta)=b\sin^{M_2}\theta\cos^{N_2}\theta$ is an even function in $C([0,2\pi])$. Thus, the result is a consequence of the Theorem \ref{P1}. Now, we observe that the system \eqref{eq.15} does not satisfy the reversibility criterion of Poincar\'e, when $N_2=0$.
\end{example}

A question on Theorem \ref{P1} is if the conditions $P_n(\cos\theta,\sin\theta)$ is an odd function and $Q_n(\cos\theta,\sin\theta)$ is an even function are also necessary conditions for (\ref{eq.13}) to have a center at the origin. The answer for this question is no, as we show in the example below. 

\begin{example}\label{ex2}
Let 
\begin{equation}\label{eq.16}
\displaystyle{\begin{array}{ccc}
\dot{x} & = &  -y+P_{n}(x,y)\\
 \dot{y} & = &  x+Q_{n}(x,y),\\
\end{array}}
 \end{equation}
where $P_n(x,y)=yP_{n-1}(x,y)$, $Q_n(x,y)=-xP_{n-1}(x,y)$ and $P_{n-1}(1,y)$ has all monomials with odd degree. We observe that $P_n(\cos\theta,\sin\theta)$ is an even function and $Q_n(\cos\theta,\sin\theta)$ is an odd function. As already discussed, the planar system (\ref{eq.16}) becomes
\begin{equation}\label{eq.17}
\displaystyle{\frac{d\gamma}{d\theta} =-(n-1)A(\theta)B(\theta)\gamma^3 + [(n-1)A(\theta)-B'(\theta)]\gamma^{2}},
\end{equation}
where
\begin{equation}\label{eq.18}
\displaystyle{\begin{array}{ccl}
A(\theta) & = &  \cos\theta\sin\theta P_{n-1}(\cos\theta,\sin\theta) -\cos\theta \sin\theta P_{n-1}(\cos\theta,\sin\theta)=0,\\
B(\theta) & = & -\cos^2\theta P_{n-1}(\cos\theta,\sin\theta) - \sin^2\theta P_{n-1}(\cos\theta,\sin\theta)\\
&=&-P_{n-1}(\cos\theta,\sin\theta).\\
\end{array}}
 \end{equation}
According to \cite{LLLZ}, Proposition 2.2 (a), the planar system (\ref{eq.16}) has a center at the origin. This shows that the conditions $P_n(\cos\theta,\sin\theta)$ is an odd function and $Q_n(\cos\theta,\sin\theta)$ is an even function are not necessary conditions for (\ref{eq.13}) to have a center at the origin.
\end{example}

In \cite{LLZ}, Llibre \textit{et al} showed a class of planar systems that have center in the origin. One of the conditions for a system to belong to this class is 
\begin{center}
\begin{equation}\label{eq.50}
f'(\theta)g(\theta)-f(\theta)g'(\theta)=ag(\theta)^3, 
\end{equation}
\end{center}
for some $a\in \mathbb{R}$. The example below shows that some planar systems satisfy the hypothesis of the Theorem \ref{P1}, but do not satisfy (\ref{eq.50}).

\begin{example}\label{ex1}
 Let the planar system 
\begin{equation}
\displaystyle{\begin{array}{ccc}
\dot{x} & = &  -y+2x^2y\\
 \dot{y}& = &  x+xy^{2}.\\
\end{array}}
 \end{equation} 
In this case, after calculations, we obtain
\[
f'(\theta)g(\theta)-f(\theta)g'(\theta)=-336\sin(\theta)\cos(\theta)^9+144\sin(\theta)\cos(\theta)^7+192\sin(\theta)\cos(\theta)^{11}
\] 
and 
\[
g(\theta)^3=(6\cos(\theta)^3\sin(\theta))^3.  
\]
It is easy to verify that $f'(\theta)g(\theta)-f(\theta)g'(\theta)\neq ag(\theta)^3$ for all $a\in\mathbb{R}$.
\end{example}

\begin{remark} 
Note that, for $A(\theta)$ in the assumptions of Theorem \ref{P1}, we have $\int_{0}^{2\pi}A(\theta)d\theta=0$. If this integral is different from zero, then the system \eqref{eq.1} has a focus at the origin. Indeed, according to \cite{Co}, the system \eqref{eq.1} is nondegenerate and quasi homogeneous. In the same paper, Conti \cite[Theorem 7.1, p. 219]{Co} proved the origin of \eqref{eq.1} is a center or a focus. 
Moreover, using the classical results of  Alwash and  Lloyd \cite{AL}, it is well known that for \eqref{eq.6}, a necessary condition to have a center is $\int_{0}^{2\pi}A(\theta)d\theta=0$.
\end{remark} 

%
%
%
%
%
\subsection{Results on Abel equations}

Consider the Abel equation
\begin{equation}\label{A1}
	x'(t)=f(t)x^3(t)+g(t)x^2(t), 
\end{equation}
$t \in [-a,a]$, where $a>0$ is a constant, $f$ and $g$ are continuous functions. We obtain conditions on $f$ and $g$ coefficients of Abel equation that ensure the existence of a center in $x=0$.

We state the main results of this section.
\begin{theorem}\label{P3}
Suppose that $f$ is an odd continuous function in $C([-a,a])$. Then, there are infinitely many closed even solutions for \eqref{A1} near the zero solution if and only if $g$ is an odd continuous function in $C([-a,a])$.  
\end{theorem}

\begin{theorem}\label{P4}
Suppose that $g$ is an odd continuous function in $C([-a,a])$. Then, there are infinitely many closed even solutions for \eqref{A1} near the zero solution if and only if $f$ is an odd continuous function in $C([-a,a])$.  
\end{theorem}

A conclusion which follows these theorems is that if $f$ and $g$ are
 odd continuous functions, each solution for \eqref{A1} near the solution $x=0$ is a closed even solution.

The next theorem proves a sufficient condition for the existence of a center to the Abel equation \eqref{A1}. 

\begin{theorem}\label{P5}
Suppose that $f$ and $g$ are odd continuous function in $C([-a,a])$. Then, the origin $x=0$ is a center of the Abel equation \eqref{A1}.  
\end{theorem}

To the particular case where $g(t)=2t$ and $f(t)$ are odd polynomials in $[-a,a]=[-1,1]$, this result is a consequence of \cite[Theorem 55, p.110]{LY}, where the authors prove the equivalence to the existence of a center. Our Theorem \ref{P5} is more general, since the assumptions on $g$ and $f$ include the polynomial case, but we prove only the sufficient conditions on $g$ and $f$ for the existence of a center. Indeed, as seen in the example below, the converse is not true.

\begin{example}\label{ex3} 
In \cite{AGG}, Alvarez, Gasull and Giacomini proved the following result. Consider the Abel equation
\begin{equation}\label{eq1}
x'=(a_0+a_1cos(2\pi t)+a_2sen(2\pi t))x^3+(b_0+b_1cos(2\pi t)+b_2sen(2\pi t))x^2
\end{equation}
where $a_0,~ a_1,~ a_2,~ b_0,~ b_1$ and $b_2$ are arbitrary real numbers. According to \cite{AGG}, for $a_0 = b_0 = a_2b_1 - a_1b_2 = 0$, the equation \eqref{eq1} has a center at $x = 0$. 
\end{example}

Considering $a_0=b_0=a_2=b_2=0$ and $a_1=b_1=1$, we obtain the following Abel equation

\begin{equation}\label{eq2}
x'=cos(2\pi t)x^3+cos(2\pi t)x^2
\end{equation}

Note that, in this case, the coefficients of Abel equation , $f(t)=cos(2\pi t)$ e $g(t)=cos(2\pi t)$, are both even functions. This shows that the reciprocal of Theorem \ref{P5} is not true. 

When $f$ and $g$ are odd polynomials, then it can be written $f(t)=t\hat{f}(t^2)$ and $g(t)=t\hat{g}(t^2)$ and the results of Theorem \ref{P5} is a consequence of the results, due to
Alwash and Lloyd \cite{AL}, which we present below.
\begin{proposition}[\cite{AL}]
Assume $f, g \in C([a,b])$ to be expressed by 
\begin{equation}\label{Composition}
	f(t)=\hat{f}(\sigma(t))\sigma'(t), \,\, g(t)=\hat{g}(\sigma(t))\sigma'(t)
\end{equation}
for some continuous functions $\hat{f}, \hat{g}$ and a continuously differentiable function $\sigma$, which is closed, i.e., $\sigma(a)=\sigma(b)$. Then, the Abel equation
\[
x'(t)=f(t)x^3(t)+g(t)x^2(t), \,\, t \in [a,b]
\]
has a center $x=0$.
\end{proposition}

In the proof of the Theorem \ref{P5} we do not use the composition condition \eqref{Composition}. The proof of the Theorem \ref{P5} follows of the Theorem \ref{P3} and Lemma \ref{LY}, when $f$ and $g$ are only odd continuous functions. Thus, we have another result with conclusions similar to those of Alwash and Lloyd.   

\section{Proof  of the theorems}

\subsection{ Preliminary results}

 We claim that a solution of \eqref{A1} is equivalent to a solution of the integral equation
\[
x(t)=\frac{\rho}{\displaystyle 1-\rho\int_{-a}^t(f(s)x(s)+g(s))ds}, t \in [-a,a]
\]
where $x(-a)=\rho$, for $\rho$ small enough, such that $\rho\int_{-a}^t(f(s)x(s)+g(s))ds <1$ for all $t\in[-a,a]$. Indeed, the equation \eqref{A1} is equivalent to
\[
-(x^{-1}(t))'=f(t)x(t)+g(t).
\]  
By integration from $-a$ to $t$ with $t \in [-a,a]$ we get
\[
\frac{1}{x(t)}-\frac{1}{x(-a)}=-\int_{-a}^t(f(s)x(s)+g(s))ds.
\]
After some computations, we obtain
\[
x(-a)=x(t)\left[ 1 - x(-a)\int_{-a}^t(f(s)x(s)+g(s))ds\right].
\]
If $x(-a)=\rho$ satisfies $\rho\int_{-a}^t(f(s)x(s)+g(s))ds <1$, for each $t \in [-a,a]$, it follows the claim.

Let $B_M(0)=\{x \in C([-a,a]); \|x\|_{\infty}\leq M\}$ be the closed ball in $C([-a,a])$, $F=\max_{t \in [-a,a]}|f(t)|$ and $G=\max_{t \in [-a,a]}|g(t)|$. Now, we define the operator $\Omega: C([-a,a]) \to C([-a,a])$ by 
\begin{equation}\label{A3}
\Omega(x)(t)=\frac{\rho}{\displaystyle 1-\rho\int_{-a}^t(f(s)x(s)+g(s))ds}, t \in [-a,a].
\end{equation}
If 
\begin{equation}\label{A4}
0\leq \rho<\frac{1}{4a(FM+G)}, 
\end{equation}
the operator \eqref{A3} is well defined. Indeed, by \eqref{A4}
\[
\rho\int_{-a}^t(f(s)x(s)+g(s))ds \leq 2a\rho(FM+G) <\frac{1}{2}.
\]
Therefore, 
\begin{equation}\label{A5}
	\frac{1}{2}< 1-\rho\int_{-a}^t(f(s)x(s)+g(s))ds
\end{equation}
and $\Omega$ is well defined.

Note that a fixed point of $\Omega$ is a solution of \eqref{A1}.

\begin{lemma}\label{L1}
The operator $\Omega$ is continuous.
\end{lemma}
\proof{} For each $x, y \in C([-a,a])$, we have
\[
\begin{array}{l}
|\Omega(x)(t) - \Omega(y)(t)|\\
=\displaystyle \left| \frac{\rho^2}{\left(1-\rho\int_{-a}^t(f(s)x(s)+g(s))ds\right)\left(1-\rho\int_{-a}^t(f(s)y(s)+g(s))ds\right)}\right| \\
\displaystyle \times\left[\int_{-a}^t(f(s)y(s)+g(s))ds - \int_{-a}^t(f(s)x(s)+g(s))ds\right]\\
\leq \displaystyle 4\rho^2\int_{-a}^a|f(s)(y(s)-x(s))|ds\\
\leq 8a\rho^2F\|y-x\|_{\infty}\\
\end{array}
\]
for each $t \in [-a,a]$. Hence
\[
\|\Omega(x) - \Omega(y)\|_{\infty}\leq  8a\rho^2F\|y-x\|_{\infty}, \forall x,y \in C([-a,a]), 
\]
and this proved the continuity.

\qed

Now, we define the restriction $\Omega_M=\Omega|_{B_M(0)}: B_M(0) \to C([-a,a])$. By Lemma \ref{L1}, $\Omega_M$ is continuous.
 
\begin{lemma}\label{L2.2}
We have that $\Omega_M: B_M(0) \to C([-a,a])$ is compact and if 
\begin{equation}\label{A7.2}
0<\rho<\min\left\{\frac{M}{2},\frac{1}{4a(FM+G)}\right\}, 
\end{equation}
we have
\[
\Omega_M(B_M(0)) \subset B_M(0).
\]
\end{lemma}
\proof{} We have that  $\Omega_M: B_M(0) \to C([-a,a])$ is bounded. Indeed, by \eqref{A5} and \eqref{A7.2}, we obtain
\[
|\Omega_M(x)(t)|\leq  2\rho, \forall x \in B_M(0), t \in [-a,a]. 
\]
Hence 
\begin{equation}\label{A6}
\|\Omega_M(x)\|_{\infty}\leq  2\rho, \forall x \in B_M(0). 
\end{equation}
For each $t, \xi \in [-a,a]$, that we can consider $t>\xi$, we have
\[
\begin{array}{l}
|\Omega_M(x)(t) - \Omega_M(x)(\xi)|\\
=\displaystyle \left| \frac{\rho^2}{\left(1-\rho\int_{-a}^t(f(s)x(s)+g(s))ds\right)\left(1-\rho\int_{-a}^{\xi}(f(s)x(s)+g(s))ds\right)}\right| \\
\displaystyle \times\left[\int_{-a}^t(f(s)x(s)+g(s))ds - \int_{-a}^{\xi}(f(s)x(s)+g(s))ds\right]\\
\leq\displaystyle 4\rho^2\int_{\xi}^t|f(s)x(s)+g(s)|ds\\
\leq 4\rho^2(FM+G)|t-\xi|, \quad \forall x \in B_M(0).\\
\end{array}
\]
Therefore, $\Omega_M(B_M(0))$ is an equicontinuous subset of $C([-a,a])$. By Ascoli-Arzela Theorem, see \cite[ p.772]{Ze}, $\Omega_M: B_M(0) \to C([-a,a])$ is compact.

Now, by \eqref{A7.2}
\[
\|\Omega_M(x)\|_{\infty}\leq  2\rho<M, \forall x \in B_M(0).
\]
Therefore, $\Omega_M: B_M(0) \to B_M(0)$ is well defined.

\qed

Now, we define the closed subspace of $C([-a,a])$ defined by
\[
E=\{x \in C([-a,a]);  \quad x \mbox{ is even}\}.
\]
Also, we define the restriction $\Omega_E=\Omega|_E: E \to C([-a,a])$. Let $B_M^E(0)=\{x \in E; \|x\|_{\infty}\leq M\}$ be the closed ball in $E$. By Lemma \ref{L1}, 
\[
\Omega_E: B_M^E(0) \to C([-a,a])
\]
is continuous.
 
\begin{lemma}\label{L2}
We have that $\Omega_E: B_M^E(0) \to C([-a,a])$ is compact and if 
\begin{equation}\label{A7}
0\leq \rho<\min\left\{\frac{M}{2},\frac{1}{4a(FM+G)}\right\}, 
\end{equation}
we have
\[
\Omega_E(B_M^E(0)) \subset B_M^E(0).
\]
\end{lemma}
\proof{} We have that  $\Omega_E: B_M^E(0) \to C([-a,a])$ is bounded. Indeed, by \eqref{A5} and \eqref{A7}, we obtain
\[
|\Omega_E(x)(t)|\leq  2\rho, \forall x \in B_M^E(0), t \in [-a,a]. 
\]
Hence 
\begin{equation}\label{A6}
\|\Omega_E(x)\|_{\infty}\leq  2\rho, \forall x \in B_M^E(0). 
\end{equation}
For each $t, \xi \in [-a,a]$, that we can consider $t>\xi$, we have
\[
\begin{array}{l}
|\Omega_E(x)(t) - \Omega_E(x)(\xi)|\\
=\displaystyle \left| \frac{\rho^2}{\left(1-\rho\int_{-a}^t(f(s)x(s)+g(s))ds\right)\left(1-\rho\int_{-a}^{\xi}(f(s)x(s)+g(s))ds\right)}\right| \\
\displaystyle \times\left[\int_{-a}^t(f(s)x(s)+g(s))ds - \int_{-a}^{\xi}(f(s)x(s)+g(s))ds\right]\\
\leq\displaystyle 4\rho^2\int_{\xi}^t|f(s)x(s)+g(s)|ds\\
\leq 4\rho^2(FM+G)|t-\xi|, \quad \forall x \in B_M^E(0).\\
\end{array}
\]
Therefore, $\Omega_E(B_M^E(0))$ is an equicontinuous subset of $C([-a,a])$. By Ascoli-Arzela Theorem, see \cite[ p.772]{Ze}, $\Omega_E: B_M^E(0) \to C([-a,a])$ is compact.

As the functions $f,g$ are odd and $x$ is even, we have that $fx+g$ is an odd function. Therefore,
\[
\int_{-a}^t(f(s)x(s)+g(s))ds
\]
is an even function. Hence, for each $x \in E$,
\[
\Omega_E(x)(t)=\frac{\rho}{\displaystyle 1-\rho\int_{-a}^t(f(s)x(s)+g(s))ds},
\]
is an even function in $[-a,a]$, that is, $\Omega_E(x) \in E$ for each $x \in E$. Now, by \eqref{A7}
\[
\|\Omega_E(x)\|_{\infty}\leq  2\rho<M, \forall x \in B_M^E(0).
\]
Therefore, $\Omega_E: B_M^E(0) \to B_M^E(0)$ is well defined.

\qed

\subsection{Proof of the Theorem \ref{P3}}

Suppose that $g$ is an odd function. Let us suppose,
\[
0\leq \rho<\min\left\{\frac{M}{2},\frac{1}{4a(FM+G)}\right\}.
\] 
It follows from Lemmas \ref{L1} and \ref{L2} that  $\Omega_E: B_M^E(0) \to B_M^E(0)$ is well defined, continuous and compact, where 
\[
\Omega_E: B_M^E(0) \to C([-a,a])
\]
and
\[
B_M^E(0)=\{x \in E; \|x\|_{\infty}\leq M\}.
\]
By the Schauder fixed point Theorem, see \cite[ p.56]{Ze}, $\Omega_E$ has a fixed point $x$, such that,
\[
\Omega_E(x)(t)=x(t)=\frac{\rho}{\displaystyle 1-\rho\int_{-a}^t(f(s)x(s)+g(s))ds}
\]
and
\[
x(-a)=\rho,
\]
for each $\rho\in \left[0,\min\left\{\frac{M}{2},\frac{1}{4a(FM+G)}\right\}\right)$. Since $\Omega_E(x)$ is an even function, we have
\[
x(-a)=x(a)
\]
and there are infinitely many closed even solutions for \eqref{A1} near the zero solution.

Now, suppose that there are infinitely many closed even solutions for \eqref{A1} near the zero solution. Note that a solution of \eqref{A1} is equivalent to the solution of the integral equation
\[
x(t)=\frac{\rho}{\displaystyle 1-\rho\int_{-a}^t(f(s)x(s)+g(s))ds}
\]
where $x(-a)=\rho$ is small enough, which is equivalent to
\[
\int_{-a}^tg(s)ds = -\frac{1}{x(t)} + \frac{1}{\rho} - \int_{-a}^tf(s)x(s)ds.
\]
Since $f$ is odd and $x(t)$ is even, we obtain $\int_{-a}^tf(s)x(s)ds$ is even. Therefore, $\int_{-a}^tg(s)ds$ is even and consequently we obtain that $g(t)$ is odd, and the theorem is proved.

 \qed

\begin{remark}
Similar arguments are true when  
\[
-\min\left\{\frac{M}{2},\frac{1}{4a(FM+G)}\right\} < \rho\leq 0.
\] 

\end{remark}

\subsection{Proofs of the Theorems \ref{P4} and \ref{P5}}

To the proof of the Theorem \ref{P5},  we apply the following result of Yang Lijun and Tang Yun \cite[Lemma 5.2, p 108]{LY}.

\begin{lemma}\label{LY}
The origin $x=0$ is a center of the Abel equation \eqref{A1}, if and only if 
\begin{equation}\label{condit1}
\int_{-a}^ag(t)dt=0 \mbox{ and } \int_{-a}^af(t)x(t,\rho)dt=0, \,\, |\rho|<\rho_0
\end{equation}
to $\rho_0$ small enough, where 
\[
x(t,\rho)=\frac{\rho}{\displaystyle 1-\rho\int_{-a}^t(f(s)x(s)+g(s))ds}.
\]
\end{lemma}  
\proof{ of Theorem \ref{P4}:} Suppose that $f$ is an odd function. The proof that there are infinitely many closed even solutions for \eqref{A1} near the zero solution is analogous to the Theorem \ref{P3}.

Now, suppose that there are infinitely many closed even solutions for \eqref{A1} near the zero solution. Note that a solution of \eqref{A1} is equivalent to the solution of the integral equation
\[
x(t)=\frac{\rho}{\displaystyle 1-\rho\int_{-a}^t(f(s)x(s)+g(s))ds}
\]
where $x(-a)=\rho$ is small enough, which is equivalent to
\[
\int_{-a}^tf(s)x(s)ds = -\frac{1}{x(t)} + \frac{1}{\rho} - \int_{-a}^tg(s)ds.
\]
Since $g$ is odd and $x(t)$ is even, we obtain $\int_{-a}^tg(s)ds$ is even. Therefore, $\int_{-a}^tf(s)x(s)ds$ is even and consequently, we obtain that $f(t)x(t)$ is odd. Since $x(t)$ is even, we conclude that $f(t)$ is odd and the theorem is proved.

 \qed

\proof{ of Theorem \ref{P5}:} It follows from Theorem \ref{P3} that there are infinitely many closed even solutions for \eqref{A1} near the zero solution. Now, we claim that, if $y(t)$ is a solution for \eqref{A1} and satisfies $y(-a)=\rho$, with $0\leq \rho < \min\left\{\frac{M}{2},\frac{1}{4a(FM+G)}\right\}$, then $y(t)$ is a closed even solution. Indeed, by Theorem \ref{P3}, there is $\bar{x}(t)$ a closed even solution
for \eqref{A1} such that $\bar{x}(-a)=\rho$. Therefore, $y(-a)=\bar{x}(-a)$ and, by uniqueness, we obtain that $y=\bar{x}$. Consequently, this proves the claim. In conclusion, if $f$ and $g$ are odd continuous functions, each solution for \eqref{A1} near the solution $x=0$ is a closed even solution.

Since
\[
\int_{-a}^ag(t)dt=0 \mbox{ and } \int_{-a}^af(t)x(t,\rho)dt=0, \,\, |\rho|<\rho_0=\min\left\{\frac{M}{2},\frac{1}{4a(FM+G)}\right\}
\]
it follows from Lemma \ref{LY} that the zero solution is a center of the Abel equation \eqref{A1}.

\qed


\

%
%
%
%

\subsection{Proofs of the Theorem \ref{P1} } 

Recall that the Abel equation associated to the planar system \eqref{eq.13} is

\begin{equation}\label{eq.30}
\displaystyle{\frac{d\gamma}{d\theta} =-(n-1)A(\theta)B(\theta)\gamma^3 + [(n-1)A(\theta)-B'(\theta)]\gamma^{2}},
\end{equation}

where
\begin{equation}\label{eq.31}
\displaystyle{\begin{array}{ccl}
A(\theta) & = &  \cos\theta P_n(\cos\theta,\sin\theta) + \sin\theta Q_n(\cos\theta,\sin\theta),\\
B(\theta) & = & \cos\theta Q_n(\cos\theta,\sin\theta) - \sin\theta P_n(\cos\theta,\sin\theta).\\
\end{array}}
 \end{equation}
As $P_n(\cos\theta,\sin\theta)$ is odd and $Q_n(\cos\theta,\sin\theta)$ is even, it follows that $f(\theta)=-(n-1)A(\theta)B(\theta)$ and $g(\theta)=(n-1)A(\theta)-B'(\theta)$ are odd continuous functions in $C([0,2\pi])$. By Theorem \ref{P5}, $\gamma=0$ is a center of the equation \eqref{eq.30}. By the equivalence stated at the introduction, the planar system \eqref{eq.13} has a center at the origin.
\qed

%
%



\end{document}